\documentclass[11pt]{article}
\usepackage{amsfonts,amsmath,amssymb,amscd,amsthm,mathrsfs,hyperref}
\title{TOPOLOGICAL AND DIFFERENTIABLE RIGIDITY OF SUBMANIFOLDS IN SPACE FORMS\footnote{
\indent Research supported by the NSFC, Grant No.
11071211; the Trans-Century Training Programme Foundation for
Talents by the Ministry of Education of China.}}
\author{HONG-WEI XU AND JUAN-RU GU}
\date{}
\numberwithin{equation}{section}
\begin{document}

\maketitle
\begin{abstract}
Let $F^{n+p}(c)$ be an $(n+p)$-dimensional simply connected space
form with nonnegative constant curvature $c$.
 We prove that if $M^n(n\geq4)$ is a compact submanifold in
$F^{n+p}(c)$, and if $Ric_M>(n-2)(c+H^2),$ where $H$ is the mean
curvature of $M$, then $M$ is homeomorphic to a sphere. We also show
that the pinching condition above is sharp. Moreover, we obtain a
new differentiable sphere theorem for submanifolds
with positive Ricci curvature.\\\\
\end{abstract}
 \section{Introduction}
 \hspace*{5mm}The investigation of curvature and topology of
Riemannian manifolds and submanifolds is one of the main stream in
global differential geometry. In 1951, Rauch first proved a
topological sphere theorem for positive pinched compact manifolds.
During the past sixty years, there are many progresses on sphere
theorems for Riemannian manifolds and submanifolds
\cite{Berger,Brendle0,Brendle1,Hamilton,Hamilton2,Shiohama2}.
 Recently B\"{o}hm and Wilking
\cite{BW} proved that every manifold with 2-positive curvature
operator must be diffeomorphic to a space form. More recently,
Brendle and Schoen \cite{Brendle2} proved the remarkable
differentiable sphere theorem for manifolds with pointwise
1/4-pinched curvatures. Moreover, Brendle and Schoen \cite{Brendle3}
obtained a differentiable rigidity theorem for compact manifolds
with weakly 1/4-pinched curvatures in the pointwise sense. The
following important convergence result for the normalized Ricci flow
in higher dimensions, initiated by Brendle and Schoen
\cite{Brendle2}, was
finally verified by Brendle \cite{Brendle}.\\\\
\textbf{Theorem A.} \emph{Let $(M,g_0)$ be a compact Riemannian
manifold of dimension $n(\geq 4)$. Assume that
$$
R_{1313}+\lambda^2R_{1414}+R_{2323}+\lambda^2R_{2424}-2\lambda
R_{1234}>0
$$
for all orthonormal four-frames $\{e_1,e_2,e_3,e_4\}$ and all
$\lambda \in [-1,1]$. Then the normalized Ricci flow with initial
metric $g_0$
$$\frac{\partial}{\partial t}g(t)=-2Ric_{g(t)}+\frac{2}{n}r_{g(t)}g(t)$$
exists for all time and converges to a constant curvature metric as
$t\rightarrow\infty$. Here $r_{g(t)}$ denotes the mean value of the
scalar curvature of $g(t)$.}\\

Let $M^{n}$ be an $n(\geq2)$-dimensional submanifold in an
$(n+p)$-dimensional Riemannian manifold $N^{n+p}$. Denote by $H$ and
$S$ the mean curvature and the squared length of the second
fundamental form of $M$, respectively. After the pioneering rigidity
theorem for minimal submanifolds in a sphere due to Simons
\cite{Simons}, Lawson \cite{Lawson} and Chern-do Carmo-Kobayashi
\cite{Chern}, Yau \cite{Yau} and Ejiri \cite{Ejiri} obtained three
important rigidity theorems for oriented compact minimal
submanifolds in $S^{n+p}$. In 1990, the first named author
\cite{Xu0} proved the following generalized Simons-Lawson-Chern-do
Carmo-Kobayashi theorem for compact
submanifolds with parallel mean curvature in a sphere.\\\\
\textbf{Theorem B.} \emph{Let $M$ be an $n$-dimensional oriented
compact submanifold with parallel mean curvature in an
$(n+p)$-dimensional unit sphere $S^{n+p}$. If $S\leq C(n,p,H),$ then
$M$ is either a totally umbilic sphere $S^{n}(\frac{1}
{\sqrt{1+H^{2}}})$, a Clifford hypersurface in an $(n+1)$-sphere, or
the Veronese surface in $S^{4}(\frac{1} {\sqrt{1+H^{2}}})$. Here the
constant $C(n,p,H)$ is defined by}
$$C(n,p,H)=\left\{\begin{array}{llll} \alpha(n,H),&\mbox{\ $for$\ } p=1, \mbox{\ $or$\ } p=2 \mbox{\  $and$\ }
 H\neq0,\\
 \frac{n}{2-\frac{1}{p}}, &\mbox{\ $for$\ } p\geq2 \mbox{\ $and$\ } H=0, \\
 \min\Big\{\alpha(n,H),\frac{n+nH^2}{2-\frac{1}{p-1}}+nH^2\Big\},&\mbox{\ $fo$r\ } p\geq 3 \mbox{\  $and$\ }
 H\neq0,
 \end{array} \right.$$
 $$\alpha(n,H)=n+ \frac{n^{3}}{2(n-1)}H^{2} -
\frac{n(n-2)}{2(n-1)}\sqrt{n^{2}H^{4}+4(n-1)H^{2}}.$$ \hspace*{5mm}
Later, the above pinching constant $C(n,p,H)$ was improved, by Li-Li
\cite{Li} for $H=0$ and by Xu \cite{Xu} for $H\neq 0$, to
$$C'(n,p,H)=\left\{\begin{array}{llll} \alpha(n,H),&\mbox{\ for\ } p=1, \mbox{\ or\ } p=2 \mbox{\ and\ }
 H\neq0,\\
\min\Big\{\alpha(n,H),\frac{1}{3}(2n+5nH^2)\Big\},&\mbox{\
otherwise.\ }
\end{array} \right.$$

Using nonexistence for stable currents on compact submanifolds of a
sphere and the generalized Poincar$\acute{e}$ conjecture in
dimension $n(\geq5)$ verified by Smale, Lawson and Simons
\cite{Lawson2} proved that if $M^{n}(n\ge 5)$ is an oriented compact
submanifold in $S^{n+p}$, and if $S<2\sqrt{n-1}$, then $M$ is
homeomorphic to a sphere. Let $F^{n+p}(c)$ be an $(n+p)$-dimensional
simply connected space form with nonnegative constant curvature $c$.
Putting
$$\alpha(n,H,c)=nc+
\frac{n^{3}}{2(n-1)}H^{2}-\frac{n(n-2)}{2(n-1)}\sqrt{n^{2}H^{4}+4(n-1)cH^{2}},$$
we have $\min_{H}\alpha(n,H,c)=2\sqrt{n-1}c.$ Motivated by the
rigidity theorem above,
Shiohama and Xu \cite{Shiohama} improved Lawson-Simons' result and proved the optimal sphere theorem.\\\\
\textbf{Theorem C.} \emph{Let $M^{n}(n\ge4)$ be an oriented complete
submanifold in $F^{n+p}(c)$ with $c\geq 0$. Suppose that
$\sup_{M}(S-\alpha(n,H,c))<0.$ Then $M$ is homeomorphic
to a sphere.}\\

The study of differentiable pinching problem for
submanifolds was initiated by Xu and Zhao \cite{XZ}. Making use of
the convergence results of Hamilton and Brendle for Ricci flow and
the Lawson-Simons formula for the nonexistence of stable currents,
Gu and Xu \cite{XG3,XG1} proved the
following differentiable sphere theorem for submanifolds in space forms.\\\\
\textbf{Theorem D.} \emph{Let $M$ be an $n(\geq4)$-dimensional
oriented complete submanifold in $F^{n+p}(c)$ with $c\ge0$. Assume
that
$S\leq 2c+\frac{n^2H^2}{n-1}$, where $c+H^2>0$. We have\\
\hspace*{2mm}$(i)$ If $c=0$, then M is either diffeomorphic to $S^n$, $\mathbb{R}^n$, or locally isometric to $S^{n-1}(r)\times \mathbb{R}$.\\
\hspace*{2mm}$(ii)$ If M is compact, then M is diffeomorphic to
$S^n$ .}\\

When $M$ is compact and $c=0$, Andrews-Baker \cite{Andrews} obtained
the same sphere theorem by the convergence result for mean curvature
flow independently.

 Recently, Xu and Gu \cite{XG2} proved
 the following generalized Ejiri rigidity theorem for compact submanifolds with parallel mean curvature in space forms.\\\\
\textbf{Theorem E.} \emph{Let $M$ be an $n(\geq3)$-dimensional
 oriented compact submanifold with parallel mean curvature in $F^{n+p}(c)$ with  $c+H^2>0$.
 If  $$Ric_{M}\geq(n-2)(c+H^2),$$
then M is either a totally umbilic sphere
$S^n(\frac{1}{\sqrt{c+H^2}})$, a Clifford hypersurface
$S^{m}\big(\frac{1}{\sqrt{2(c+H^2)}}\big)\times
S^{m}\big(\frac{1}{\sqrt{2(c+H^2)}}\big)$ in the totally umbilic
sphere $S^{n+1}(\frac{1}{\sqrt{c+H^2}})$ with $n=2m$, or
$\mathbb{C}P^{2}(\frac{4}{3}(c+H^2))$ in
$S^7(\frac{1}{\sqrt{c+H^2}})$. Here
$\mathbb{C}P^{2}(\frac{4}{3}(c+H^2))$ denotes the $2$-dimensional
complex projective space minimally immersed in
$S^7(\frac{1}{\sqrt{c+H^2}})$ with constant holomorphic sectional
curvature $\frac{4}{3}(c+H^2)$.}\\

The purposes of the present paper is to investigate rigidity of
topological and differentiable structures of compact submanifolds.
Our paper is organized as follows. Some notations and lemmas are
prepared in Section 2.  In Section 3, we prove that if $M$ is an
$n(\geq4)$-dimensional compact submanifold in $F^{n+p}(c)$ with
$c\geq 0$, and
 if  $Ric_M>(n-2)(c+H^2),$
then $M$ is homeomorphic to a sphere. We then give an example to
show that the pinching condition above is sharp. Moreover, we obtain
a new differentiable sphere theorem for compact submanifolds with
positive
Ricci curvature in a space form.\\\\

\section{Notations and lemmas}\hspace*{5mm}Throughout this paper let $M^n$ be an $n$-dimensional compact
submanifold in an $(n+p)$-dimensional Riemannian manifold $N^{n+p}$.
We shall make use of the following convention on the range of
indices:
$$1\leq A,B,C,\cdots\leq n+p,\ 1\leq i,j,k,\cdots\leq n,\
  n+1\leq \alpha,\beta,\gamma,\cdots\leq n+p.$$
For an arbitrary fixed point $x\in M\subset N$, we choose an
orthonormal local frame field $\{e_{A}\}$ in $N^{n+p}$ such that
$e_{i}$'s are tangent to $M$. Denote by $\{\omega_A\}$ the dual
frame field of $\{e_A\}$. Let $Rm$, $h$ and $\xi$ be the Riemannian
curvature tensor, second fundamental form and mean curvature vector
of $M$ respectively, and $\overline{Rm}$ the Riemannian curvature
tensor of $N$. Then
\begin{eqnarray}&&Rm=\sum_{i,j,k,l}R_{ijkl}\omega_i\otimes\omega_j
                    \otimes\omega_k\otimes\omega_l,\nonumber\\
&&\overline{Rm}=\sum_{A,B,C,D}\overline{R}_{ABCD}\omega_A\otimes\omega_B
                         \otimes\omega_C\otimes\omega_D,\nonumber\\
&&h=\sum_{\alpha,i,j}h^{\alpha}_{ij}\omega_i\otimes
                     \omega_j\otimes e_\alpha, \,\,\xi=\frac{1}{n}\sum_{\alpha,i}h^{\alpha}_{ii}e_{\alpha},\nonumber\\
&&R_{ijkl}=\overline{R}_{ijkl}+\sum_{\alpha}\big(h_{ik}^{\alpha}h_{jl}^{\alpha}
                         -h_{il}^{\alpha}h_{jk}^{\alpha}\big),\label{Gauss}\\
&&R_{\alpha\beta kl}=\overline{R}_{\alpha\beta
kl}+\sum_{i}(h^{\alpha}_{ik}h^{\beta}_{il}-h^{\alpha}_{il}h^{\beta}_{ik}).\label{Codazzi}
\end{eqnarray}
We define
$$S=|h|^{2}, \ H=|\xi|, \ H_{\alpha}=(h^{\alpha}_{ij})_{n\times
n}.$$
Denote by $Ric(u)$ the Ricci curvature of $M$ in direction of $u\in
UM$. From the Gauss equation, we have
\begin{equation}Ric(e_i)=\sum_{j}\overline{R}_{ijij}
+\sum_{\alpha,j}\big[h_{ii}^{\alpha}h_{jj}^{\alpha}
                          -(h_{ij}^{\alpha})^2\big].\label{Ricci}\end{equation}
Set $Ric_{\min}(x)=\min_{u\in U_{x}M}Ric(u)$.
Denote by $K(\pi)$ the sectional curvature of $M$ for tangent
2-plane $\pi(\subset T_xM)$ at point $x\in M$, $\overline{K}(\pi)$
the sectional curvature of $N$ for tangent 2-plane $\pi(\subset
T_xN)$ at point $x\in N$. Set $\overline{K}_{\min}:=\min_{\pi\subset
T_{x}N}\overline{K}(\pi)$, $\overline{K}_{\max}:=\max_{\pi\subset
T_{x}N}\overline{K}(\pi)$. Then by Berger's inequality, we have
\begin{equation}|\overline{R}_{ABCD}|\leq\frac{2}{3}(\overline{K}_{\max}-\overline{K}_{\min})
\label{Berger}\end{equation}
 for all distinct indices $A$, $B$, $C$, $D$.

When the ambient manifold $N^{n+p}$ is the complete and
simply connected space form $F^{n+p}(c)$ with constant curvature $c$, we have
\begin{equation}
\overline{R}_{ABCD}=c(\delta_{AC}\delta_{BD}-\delta_{AD}\delta_{BC}).\label{Sectional}
\end{equation}
Then the scalar curvature of $M$ is given by
\begin{equation}R=n(n-1)c+n^{2}H^{2}-S.\label{Gauss2}\end{equation}

The nonexistence theorem for stable currents in a compact Riemannian
manifold $M$ isometrically immersed into $F^{n+p}(c)$ is employed to
eliminate the homology groups $H_{q}(M;\mathbb{Z})$ for $0<q<n$,
which was initiated by
Lawson-Simons \cite{Lawson2} and extended by Xin \cite{Xin}.\\\\
\textbf{Theorem 2.1.} \emph{Let $M^{n}$ be a compact submanifold in
$F^{n+p}(c)$ with $c\geq 0$. Assume that$$
\sum_{k=q+1}^{n}\sum_{i=1}^{q}[2|h(e_{i},e_{k})|^{2}-\langle
h(e_{i},e_{i}),h(e_{k},e_{k})\rangle]<q(n-q)c$$ holds for any
orthonormal basis $\{e_{i}\}$ of $T_xM$ at any point $x\in M$, where
q is an integer satisfying $0<q<n$. Then there does not exist any
stable q-currents. Moreover,
$H_{q}(M;\mathbb{Z})=H_{n-q}(M;\mathbb{Z})=0,$ and $\pi_{1}(M)=0$
when $q=1$. Here $H_{i}(M;\mathbb{Z})$ is the $i$-th homology group
of M with integer
coefficients.}\\

To prove the sphere theorems for submanifolds, we need to eliminate
the fundamental group
$\pi_1(M)$ under the Ricci curvature pinching condition, and get the following lemma.\\\\
  \textbf{Lemma 2.1.} \emph{Let $M$ be an $n(\geq4)$-dimensional
 compact submanifold in $F^{n+p}(c)$ with  $c\geq 0$.
 If the Ricci curvature of M satisfies
  $$Ric_M>\frac{n(n-1)}{n+2}(c+H^2),$$ then $H_{1}(M;\mathbb{Z})=H_{n-1}(M;\mathbb{Z})=0,$ and $\pi_1(M)=0$.}\\\\
\textbf{Proof.} From (\ref{Gauss2}) and the assumption, we have
$$S-nH^2< \frac{2n(n-1)}{n+2}(c+H^2).$$This together with (\ref{Ricci})
implies that
\begin{eqnarray} &&\sum_{k=2}^{n}[2|h(e_{1},e_{k})|^{2}-\langle
h(e_{1},e_{1}),h(e_{k},e_{k})\rangle]\nonumber\\&=&2\sum_{\alpha}\sum_{k=2}^{n}(h^{\alpha}_{1k})^{2}-
\sum_{\alpha}\sum_{k=2}^{n}h^{\alpha}_{11}h^{\alpha}_{kk}\nonumber\\
&=&\sum_{\alpha}\sum_{k=2}^{n}(h^{\alpha}_{1k})^{2}-Ric(e_1)+(n-1)c\nonumber\\
&\leq&\frac{1}{2}(S-nH^2)-Ric(e_1)+(n-1)c\nonumber\\
&<&\frac{n(n-1)}{n+2}(c+H^2)-\frac{n(n-1)}{n+2}(c+H^2)+(n-1)c\nonumber\\
&=&(n-1)c.
\end{eqnarray}
This together with Theorem 2.1 implies that
$H_{1}(M;\mathbb{Z})=H_{n-1}(M;\mathbb{Z})=0,$ and  $\pi_1(M)=0$.
This proves Lemma 2.1. \\\\

\section{Sphere theorems for submanifolds}
 \hspace*{5mm}In this section, we investigate rigidity of topological and differentiable structures of
 compact submanifolds in space forms. Motivated by Theorem E, we first prove the following
 topological sphere theorem for compact submanifolds in space forms. \\\\
\textbf{Theorem 3.1.} \emph{Let $M$ be an $n(\geq4)$-dimensional
compact submanifold in $F^{n+p}(c)$ with  $c\geq0$. If
$$Ric_M>(n-2)(c+H^2),$$  then M
is homeomorphic to a sphere.}\\\\
\textbf{Proof.} Assume that $ 2\leq q\leq \frac{n}{2}$. Setting
$$T_{\alpha}:=\frac{trH_{\alpha}}{n},$$
  we have $\sum_{\alpha}T_{\alpha}^{2}=H^2,$  and
  \begin{equation}Ric(e_i)=(n-1)c+\sum_{\alpha}\Big[nT_{\alpha}h^{\alpha}_{ii}-(h^{\alpha}_{ii})^2
  -\sum_{j\neq i}(h^{\alpha}_{ij})^2\Big].
  \end{equation}
Then we get\begin{eqnarray}
&&\sum_{k=q+1}^{n}\sum_{i=1}^{q}[2|h(e_{i},e_{k})|^{2}-\langle
h(e_{i},e_{i}),h(e_{k},e_{k})\rangle]\nonumber\\
&=&2\sum_{\alpha}\sum_{k=q+1}^{n}\sum_{i=1}^{q}(h^{\alpha}_{ik})^{2}-
\sum_{\alpha}\sum_{k=q+1}^{n}\sum_{i=1}^{q}h^{\alpha}_{ii}h^{\alpha}_{kk}\nonumber\\
&=&\sum_{\alpha}\Big[2\sum_{k=q+1}^{n}\sum_{i=1}^{q}(h^{\alpha}_{ik})^{2}-
\Big(\sum_{i=1}^{q}h^{\alpha}_{ii}\Big)\Big(trH_{\alpha}-\sum_{i=1}^{q}h^{\alpha}_{ii}\Big)\Big]\nonumber\\
&\leq&\sum_{\alpha}\Big[2\sum_{k=q+1}^{n}\sum_{i=1}^{q}(h^{\alpha}_{ik})^{2}-
nT_{\alpha}\sum_{i=1}^{q}h^{\alpha}_{ii}+q\sum_{i=1}^{q}(h^{\alpha}_{ii})^2\Big]\nonumber\\
&\leq& q\sum_{i=1}^{q}[(n-1)c-Ric(e_i)]
+n(q-1)\sum_{\alpha}\sum_{i=1}^{q}T_{\alpha}h^{\alpha}_{ii}\nonumber\\
&\leq&q^2[(n-1)(c+H^2)-Ric_{\min}]\nonumber\\
&&
-q(n-q)H^2+n(q-1)\sum_{\alpha}\sum_{i=1}^{q}T_{\alpha}(h^{\alpha}_{ii}-T_{\alpha})\nonumber\\
&\leq&q(n-q)[(n-1)(c+H^2)-Ric_{\min}]\nonumber\\
&&
-q(n-q)H^2+n(q-1)\sum_{\alpha}\sum_{i=1}^{q}T_{\alpha}(h^{\alpha}_{ii}-T_{\alpha}).\label{Estimate2}
\end{eqnarray}
 On the other hand, we obtain
\begin{eqnarray}
&&\sum_{k=q+1}^{n}\sum_{i=1}^{q}[2|h(e_{i},e_{k})|^{2}-\langle
h(e_{i},e_{i}),h(e_{k},e_{k})\rangle]\nonumber\\
&=&\sum_{\alpha}\Big[2\sum_{k=q+1}^{n}\sum_{i=1}^{q}(h^{\alpha}_{ik})^{2}-
\frac{n-q}{n}\Big(\sum_{i=1}^{q}h^{\alpha}_{ii}\Big)\Big(trH_{\alpha}-\sum_{i=1}^{q}h^{\alpha}_{ii}\Big)
\nonumber\\
&&-\frac{q}{n}\Big(\sum_{k=q+1}^{n}h^{\alpha}_{kk}\Big)\Big(trH_{\alpha}-\sum_{k=q+1}^{n}h^{\alpha}_{kk}\Big)\Big]\nonumber\end{eqnarray}
\begin{eqnarray}
&\leq&\sum_{\alpha}\Big[2\sum_{k=q+1}^{n}\sum_{i=1}^{q}(h^{\alpha}_{ik})^{2}-
(n-q)T_{\alpha}\sum_{i=1}^{q}h^{\alpha}_{ii}+\frac{q(n-q)}{n}\sum_{i=1}^{q}(h^{\alpha}_{ii})^2 \nonumber\\
&&-qT_{\alpha}\sum_{k=q+1}^{n}h^{\alpha}_{kk}+\frac{q(n-q)}{n}\sum_{k=q+1}^{n}(h^{\alpha}_{kk})^2\Big]\nonumber\\
&\leq& \frac{q(n-q)}{n}S
-\sum_{\alpha}\Big[qnT^2_{\alpha}+(n-2q)T_{\alpha}\sum_{i=1}^{q}h^{\alpha}_{ii}\Big]\nonumber\\
&\leq&q(n-q)[(n-1)(c+H^2)-Ric_{\min}]\nonumber\\
&&-q(n-q)H^2
-(n-2q)\sum_{\alpha}\sum_{i=1}^{q}T_{\alpha}(h^{\alpha}_{ii}-T_{\alpha}).\label{Estimate3}
\end{eqnarray}
 It follows from (\ref{Estimate2}), (\ref{Estimate3}) and the assumption that
\begin{eqnarray}
&&\sum_{k=q+1}^{n}\sum_{i=1}^{q}[2|h(e_{i},e_{k})|^{2}-\langle
h(e_{i},e_{i}),h(e_{k},e_{k})\rangle]\nonumber\\
 &\leq&\frac{n-2q}{q(n-2)}\Big\{q(n-q)[(n-1)(c+H^2)-Ric_{\min}]\nonumber\\
&&-q(n-q)H^2+n(q-1)\sum_{\alpha}\sum_{i=1}^{q}T_{\alpha}(h^{\alpha}_{ii}-T_{\alpha})\Big\}\nonumber\\
&&+\frac{n(q-1)}{q(n-2)}\Big\{q(n-q)[(n-1)(c+H^2)-Ric_{\min}]\nonumber\\
&&-q(n-q)H^2
-(n-2q)\sum_{\alpha}\sum_{i=1}^{q}T_{\alpha}(h^{\alpha}_{ii}-T_{\alpha})\Big\}\nonumber\\
&=&q(n-q)[(n-1)(c+H^2)-Ric_{\min}]-q(n-q)H^2\nonumber\\
&<&q(n-q)c.
\end{eqnarray}
This together with Theorem 2.1 implies that
   $$H_{q}(M;\mathbb{Z})=H_{n-q}(M;\mathbb{Z})=0,$$ for all
 $2\leq q\leq \frac{n}{2}.$

 Since $(n-2)(c+H^2)\geq \frac{n(n-1)}{n+2}(c+H^2),$ we get from the assumption and Lemma 2.1 that
 $$H_{1}(M;\mathbb{Z})=H_{n-1}(M;\mathbb{Z})=0,$$ and $M$ is simply connected.

 From above discussion, we know that $M$ is a homotopy sphere. This together with the generalized Paincar\'{e}
 conjecture implies that $M$ is a topological sphere. This completes the proof of Theorem
 3.1.\\

Theorem 3.1 improves the sphere theorem due to Vlachos
\cite{Vlachos} and Hu-Zhai \cite{Hu}. The following example shows
that our pinching condition in
Theorem 3.1 is sharp.\\\\
Example 3.1. (i) Let $M:=\mathbb{C}P^{2}(\frac{4}{3}(c+H^2))\subset
S^7(\frac{1}{\sqrt{c+H^2}})$, where $H$ is a nonnegative constant.
Then $M$ is a compact minimal submanifold in
$S^7(\frac{1}{\sqrt{c+H^2}})(\subset F^{4+p}(c))$ with $Ric_M\equiv
2(c+H^2)$. Hence $M$ is a compact submanifold in $F^{4+p}(c)$ with
constant mean curvature $H$. It is not a topological sphere.

(ii) Let $M:=S^{m}\big(\frac{1}{\sqrt{2(c+H^2)}}\big)\times
S^{m}\big(\frac{1}{\sqrt{2(c+H^2)}}\big)\subset
S^{n+1}(\frac{1}{\sqrt{c+H^2}})$ with $n=2m\geq4$, where $H$ is a
nonnegative constant. Then $M$ is a compact minimal submanifold in
$S^{n+1}(\frac{1}{\sqrt{c+H^2}})(\subset F^{n+p}(c))$ with
$Ric_M\equiv (n-2)(c+H^2)$. Hence $M$ is a compact submanifold in
$F^{n+p}(c)$ with constant mean curvature $H$. It is not a
topological
sphere.\\

In the next, we investigate differentiable pinching problem on compact
submanifolds in a Riemannian manifold, and obtain the following theorem.\\\\
 \textbf{Theorem 3.2.} \emph{Let $(M, g_0)$ be an $n(\geq 4)$-dimensional
compact submanifold in an $(n+p)$-dimensional Riemannian manifold
$N^{n+p}$. If the Ricci curvature of M satisfies $$Ric_{M}>\Big[\frac{3n^2-9n+8}{3(n-2)}\overline{K}_{\max}
                -\frac{8}{3(n-2)}\overline{K}_{\min}\Big]+\frac{n(n-3)}{n-2}H^2,$$
 then the normalized Ricci flow with initial metric $g_0$
$$\frac{\partial}{\partial t}g(t) = -2Ric_{g(t)} +\frac{2}{n}
r_{g(t)}g(t),$$ exists for all time and converges to a constant
curvature metric as $t\rightarrow\infty$. Moreover, $M$ is
diffeomorphic to a spherical space form. In particular, if M is
simply connected, then M is diffeomorphic to $S^n$.}\\\\
\textbf{Proof.} Set $T_{\alpha}= \frac{1}{n}trH_{\alpha}$. Then
$\sum_{\alpha}T_{\alpha}^2=H^2$, and
\begin{eqnarray}h_{ii}^{\alpha}h_{jj}^{\alpha}
                    &=&\frac{1}{2}[(h^{\alpha}_{ii}+h^{\alpha}_{jj}-2T_{\alpha})^2
                    -(h^{\alpha}_{ii}-T_{\alpha})^2
-(h^{\alpha}_{jj}-T_{\alpha})^2]\nonumber\\
&&+T_{\alpha}(h^{\alpha}_{ii}-T_{\alpha})+T_{\alpha}(h^{\alpha}_{jj}-T_{\alpha})+T_{\alpha}^2.\label{Estimate4}\end{eqnarray}
                    We rewrite (\ref{Ricci}) as
\begin{eqnarray}
Ric(e_i)&=&\sum_j\overline{R}_{ijij}+(n-1)H^2+(n-2)\sum_{\alpha}T_{\alpha}(h^{\alpha}_{ii}-T_{\alpha})\nonumber\\
&& -\sum_{\alpha}(h^{\alpha}_{ii}-T_{\alpha})^2-\sum_{\alpha,j\neq
i}(h^{\alpha}_{ij})^2.
\end{eqnarray}
This implies that
\begin{eqnarray}
-\sum_{\alpha}(h^{\alpha}_{ii}-T_{\alpha})^2&\geq&
Ric_{\min}-(n-1)(\overline{K}_{\max}+H^2)\nonumber\\
&&-(n-2)\sum_{\alpha}T_{\alpha}(h^{\alpha}_{ii}-T_{\alpha})
+\sum_{\alpha,j\neq i}(h^{\alpha}_{ij})^2,\label{Ricci3}
\end{eqnarray}
and
\begin{equation}
 \sum_{\alpha}T_{\alpha}(h^{\alpha}_{ii}-T_{\alpha})\geq\frac{1}{n-2}[Ric_{\min}-(n-1)(\overline{K}_{\max}+H^2)].\label{Ricci4}
\end{equation}
Suppose $\{e_1,e_2,e_3,e_4\}$ is an orthonormal
four-frame and $\lambda \in \mathbb{R}$.

From (\ref{Gauss}), (\ref{Berger}), (\ref{Estimate4}), (\ref{Ricci3}) and (\ref{Ricci4}), we
have
\begin{eqnarray}
 && R_{1313}+R_{2323}-|R_{1234}|\nonumber\\
&=&
\overline{R}_{1313}+\overline{R}_{2323}+
   \sum_{\alpha}\Big[h_{11}^{\alpha}h_{33}^{\alpha}-(h_{13}^{\alpha})^2
                    +h_{22}^{\alpha}h_{33}^{\alpha}-(h_{23}^{\alpha})^2\Big]\nonumber\\
                    &&-|\overline{R}_{1234}+\sum_{\alpha}
                    (h^{\alpha}_{13}h^{\alpha}_{24}-h^{\alpha}_{14}h^{\alpha}_{23})|\nonumber\\
  &\geq&2\overline{K}_{\min}
  -\frac{2}{3}(\overline{K}_{\max}-\overline{K}_{\min})
  -\frac{1}{2}\sum_{\alpha}\Big[3(h_{13}^{\alpha})^2
  +3(h_{23}^{\alpha})^2+(h_{14}^{\alpha})^2+(h_{24}^{\alpha})^2\Big]\nonumber\\
 &&+\sum_{\alpha}\Big[-\frac{(h^{\alpha}_{11}-T_{\alpha})^2}{2}
                    -\frac{(h^{\alpha}_{33}-T_{\alpha})^2}{2}+T_{\alpha}(h^{\alpha}_{11}-T_{\alpha})
                    +T_{\alpha}(h^{\alpha}_{33}-T_{\alpha})+T_{\alpha}^2\Big]\nonumber\\
 &&+\sum_{\alpha}\Big[-\frac{(h^{\alpha}_{22}-T_{\alpha})^2}{2}
                    -\frac{(h^{\alpha}_{33}-T_{\alpha})^2}{2}+T_{\alpha}(h^{\alpha}_{22}-T_{\alpha})
                    +T_{\alpha}(h^{\alpha}_{33}-T_{\alpha})+T_{\alpha}^2\Big]\nonumber\\
 &\geq
 &\frac{8}{3}\Big(\overline{K}_{\min}
  -\frac{1}{4}\overline{K}_{\max}\Big)-\frac{1}{2}\sum_{\alpha}\Big[3(h_{13}^{\alpha})^2
  +3(h_{23}^{\alpha})^2+(h_{14}^{\alpha})^2+(h_{24}^{\alpha})^2\Big]\nonumber\\
 &&+2[Ric_{\min}-(n-1)(\overline{K}_{\max}+H^2)]+2H^2\nonumber\\
 && +\frac{1}{2}\sum_{\alpha,j\neq 1}(h^{\alpha}_{1j})^2+\frac{1}{2}\sum_{\alpha,j\neq 2}(h^{\alpha}_{2j})^2+\sum_{\alpha,j\neq 3}(h^{\alpha}_{3j})^2
 \nonumber\\&&
 +\frac{n-4}{2}\sum_{\alpha,i\neq 1,3}T_{\alpha}(h^{\alpha}_{ii}-T_{\alpha})+\frac{n-4}{2}\sum_{\alpha,i\neq 2,3}T_{\alpha}(h^{\alpha}_{ii}-T_{\alpha})
 \nonumber\\
  &\geq
 &\frac{8}{3}\Big(\overline{K}_{\min}
  -\frac{1}{4}\overline{K}_{\max}\Big)+2H^2+(n-2)[Ric_{\min}-(n-1)(\overline{K}_{\max}+H^2)].\label{Estimate5}
 \end{eqnarray}
 Same argument implies that
 \begin{eqnarray}
 && R_{1414}+R_{2424}-|R_{1234}|\nonumber\\
  &\geq
 &\frac{8}{3}\Big(\overline{K}_{\min}
  -\frac{1}{4}\overline{K}_{\max}\Big)+2H^2+(n-2)[Ric_{\min}-(n-1)(\overline{K}_{\max}+H^2)].
 \end{eqnarray}
This together with (\ref{Estimate5}) and the assumption implies
\begin{eqnarray}&&R_{1313}+\lambda^2 R_{1414} + R_{2323} +
\lambda^2R_{2424}- 2\lambda R_{1234}\nonumber\\
&\geq&R_{1313}+ R_{2323}-|R_{1234}|+\lambda^2(R_{1414}+ R_{2424}-|R_{1234}|)\nonumber\\
&>&0.\end{eqnarray}
It follows from Theorem A that $M$ is diffeomorphic to a spherical
space form. In particular, if $M$ is simply connected, then $M$ is
diffeomorphic
to $S^n$. This completes the proof of Theorem 3.2.\\\\
 \textbf{Theorem 3.3.} \emph{Let $M$ be an $n(\geq 4)$-dimensional
compact submanifold in $F^{n+p}(c)$ with $c\ge0$.  If
  $$Ric_M >(n-2)(1+\varepsilon_n)(c+H^2),$$
 then $M$ is diffeomorphic to $S^n$. Here$$\varepsilon_n=\left\{\begin{array}{llll} 0, &\mbox{\ for\ } 4\leq n\leq 6,\\
 \frac{n-4}{(n-2)^2}, &\mbox{\
for\ }n\geq 7.
\end{array} \right.$$}\\
\textbf{Proof.} When $n=5,6$, it is well known that there is only
one differentiable structure on $S^n$. This together with Theorem
3.1 implies $M$ is diffeomorphic to $S^n$. When $n\neq 5,6$, it
follows from Theorem 3.2 that $M$ is diffeomorphic to a spherical
space form. On the other hand,  it follows from Lemma 2.1 that $M$
is simply connected. Therefore, $M$ is diffeomorphic to $S^n$. This
completes the proof of Theorem 3.3.\\\\
Remark 3.1. When $4\le n\le6$, the pinching condition in Theorem 3.3
is sharp.  When $n\ge7$, we have $0\leq \varepsilon_n<\frac{1}{n}$
and $\lim_{n\rightarrow\infty} \varepsilon_n =0.$ Therefore,
the pinching condition in Theorem 3.3 is close to the best possible.\\

Motivated by Theorem E and the sphere theorems above, we would like to
propose
the following conjecture.\\\\
\textbf{Conjecture.} \emph{Let $M$ be an $n(\geq3)$-dimensional
compact oriented submanifold in the space form $F^{n+p}(c)$ with
$c+H^2>0$.
 If  $$Ric_{M}\geq(n-2)(c+H^2),$$
 then M is diffeomorphic to either the standard
$n$-sphere $S^n$, the Clifford hypersurface
$S^{m}\big(\frac{1}{\sqrt{2}}\big)\times
S^{m}\big(\frac{1}{\sqrt{2}}\big)$ in $S^{n+1}$ with $n=2m$, or
$\mathbb{C}P^{2}$. In particular, if  $Ric_{M}>(n-2)(c+H^2)$, then $M$ is diffeomorphic to $S^n$.}\\

Theorems 3.2 and 3.3 provide partial affirmative answer to the
Conjecture.\\\\

Hong-Wei Xu

Center of Mathematical Sciences\

Zhejiang University\

Hangzhou 310027\

China

E-mail address: xuhw@cms.zju.edu.cn\\\\

Juan-Ru Gu

Center of Mathematical Sciences\

Zhejiang University\

Hangzhou 310027\

China

E-mail address: gujr@cms.zju.edu.cn


\begin{thebibliography}{bb}
\bibitem{Andrews} B. Andrews and C. Baker, Mean curvature flow of pinched submanifolds
to spheres, \emph{J. Differential Geom.}, {\bf85}(2010), 357-396.
\bibitem{Berger}M. Berger, Riemannian geometry during the second half of the twentieth century, University Lecture Series, Vol.{\bf 17},
American Mathematical Society, Providence, RI, 2000.
\bibitem{BW}C. B\"{o}hm and B.  Wilking, Manifolds with positive curvature operator are space forms,
\emph{Ann. of Math.}, {\bf167}(2008), 1079-1097.
\bibitem{Brendle}S. Brendle, A general convergence result for the Ricci flow in higher dimensions,
\emph{Duke Math. J.}, {\bf145}(2008), 585-601.
\bibitem{Brendle0}S. Brendle, Ricci Flow and the Sphere Theorem, Graduate Studies in Mathematics, Vol.{\bf111}, Americam Mathematical Society, 2010.
\bibitem{Brendle2} S. Brendle and R. Schoen, Manifolds with $1/4$-pinched curvature are space forms,
\emph{J. Amer. Math. Soc.}, {\bf22}(2009), 287-307.
\bibitem{Brendle3}S. Brendle and R. Schoen, Classification of manifolds with weakly
$1/4$-pinched curvatures, \emph{Acta Math.}, {\bf200}(2008), 1-13.
\bibitem{Brendle1}S. Brendle and R. Schoen, Sphere theorems in geometry, Surveys in Differential Geometry,
Vol.{\bf13}, 2009, 49-84.
\bibitem{Chern}S. S. Chern, M. do Carmo and S. Kobayashi, Minimal submanifolds of a sphere with second
fundamental form of constant length, in Functional Analysis and
Related Fields, Springer-Verlag, New York(1970).
\bibitem{Ejiri}N. Ejiri,  Compact minimal submanifolds of a sphere with
positive Ricci curvature, \emph{J. Math. Soc.Japan}, {\bf31}(1979), 251-256.
\bibitem{XG3} J. R. Gu and H. W. Xu, The sphere theorem for
manifolds with positive scalar curvatue, arXiv:math.DG/1102.2424.
\bibitem{Hamilton}R. Hamilton, Three manifolds with positive Ricci
curvature, \emph{J. Differential Geom.}, {\bf17}(1982), 255-306.
\bibitem{Hamilton2}R. Hamilton, Four-manifolds with positive curvature operator, \emph{J. Differential Geom.}, {\bf24}(1986), 153-179.
\bibitem{Hu}Z. J. Hu and J. S. Zhai, A sphere theorem for even-dimensional
submanifolds of the unit sphere, \emph{J. Zhenqzhou Univ. $($Nat.
Sci. Ed.$)$}, {\bf37}(2005), 1-4.
\bibitem{Lawson}B. Lawson, Local rigidity theorems for minimal hyperfaces, \emph{Ann. of Math.}, {\bf 89}(1969), 187-197.
\bibitem{Lawson2}B. Lawson and J. Simons, On stable currents and their
application to global problems in real and complex geometry,
\emph{Ann. of Math.}, {\bf98}(1973), 427-450.
\bibitem{Li}A. M. Li and J. M. Li, An intrinsic rigidity theorem for minimal
submanifolds in a sphere, \emph{Arch. Math.}, {\bf58}(1992),
582-594.
\bibitem{Petersen1}P. Petersen and T. Tao, Classification of almost
quarter-pinched manifolds, \emph{Proc. Amer. Math. Soc.},
{\bf137}(2009), 2437-2440.
\bibitem{Shiohama2}K. Shiohama, Sphere theorems, Handbook of Differential Geometry, Vol. {\bf1}, F.J.E. Dillen and
L.C.A. Verstraelen (eds.), Elsevier Science B.V., Amsterdam, 2000.
\bibitem {Shiohama} K. Shiohama and H. W. Xu, The topological sphere theorem for complete submanifolds,
\emph{Compositio Math.}, {\bf107}(1997), 221-232.
\bibitem{Simons}J. Simons, Minimal varieties in Riemannian manifolds, \emph{Ann. Math.}, {\bf 88}(1968), 62-105.
\bibitem {Sjerve}D. Sjerve, Homology spheres which are covered by
spheres. \emph{J. London Math. Soc.}, {\bf6}(1973), 333-336.
\bibitem{Vlachos}T. Vlachos, A sphere theorem for
odd-dimensional submanifolds of spheres, \emph{Proc. Amer. Math.
Soc.}, {\bf130}(2002), 167-173.
\bibitem{Xin}Y. L. Xin, Application of integral currents to
vanishing theorems, \emph{Scient. Sinica}(A), {\bf27}(1984),
233-241.
\bibitem{Xu0}H. W. Xu, Pinching theorems, global
pinching theorems, and eigenvalues for Riemannian submanifolds,
\emph{Ph.D. dissertation, Fudan University}, 1990.
\bibitem{Xu}H. W. Xu, A rigidity theorem for submanifolds with parallel mean curvature in
a sphere, \emph{Arch. Math.}, {\bf61}(1993), 489-496.
\bibitem{XG1} H. W. Xu and J. R. Gu,  An optimal differentiable sphere theorem for complete
manifolds, \emph{Math. Res. Lett.}, {\bf17}(2010), 1111-1124.
\bibitem{XG2} H. W. Xu and J. R. Gu, Rigidity of submanfolds with parallel mean curvature in space forms,
arXiv:math.DG/1105.2920v1.
\bibitem{XT} H. W. Xu and L. Tian, A differentiable sphere theorem inspired by rigidity of minimal submanifolds, \emph{to appear in Pacific J. Math.}
\bibitem{XZ}H. W. Xu and E. T. Zhao, Topological and
differentiable sphere theorems for complete submanifolds,  \emph{Comm.
Anal. Geom.}, {\bf 17}(2009), 565-585.
\bibitem{Yau}S. T. Yau, Submanifolds with constant mean curvature I, II, \emph{Amer. J. Math.}, {\bf96, 97}(1974, 1975), 346-366, 76-100.\\


\end{thebibliography}
\end{document}